\input amstex.tex
\input amsppt.sty   
\magnification = \magstep1
\parskip=4pt plus 2pt minus 2pt
\baselineskip=14pt
\parindent=6pt

\overfullrule=0pt

\NoRunningHeads
\NoBlackBoxes
\widestnumber\key{AMBDCF}
\redefine \Bbb {\bold}
\def \C {\Bbb C}
\def \N {\Bbb N}

\def \R {\Bbb R}
\def \Z {\Bbb Z}



\def \Tr{\operatorname{Tr}}

\def \Sp{\operatorname{Sp}}
\def \supp{\operatorname{supp}}

\def \Sp2{\text{Sp}(2,\R)}
\def \O2{\text{O}(2,\R)}

\def \e-phi{e^{-\phi}}



\def \ssk {\vskip 9pt}

\vsize=8.2 true in
\hsize=6.2 true in
\document
\topmatter
\title Anderson Localization\\
 for Time Periodic \\
Random Sch\"odinger Operators\endtitle
\author
Avy Soffer\\ Wei-Min Wang
\endauthor
\affil
Department of Mathematics\\
Rutgers University\\
soffer\@math.rutgers.edu\\
{}\\
UMR 8628 du CNRS \\
and\\
Department of Mathematics and Physics\\
Princeton University\\  
wmwang\@feynman.princeton.edu
\endaffil
\address New Brunswick, NJ 08903\endaddress
\address Princeton, NJ 08544\endaddress
\email soffer@math.rutgers.edu\endemail
\email wmwang@feynman.princeton.edu\endemail
\thanks  We thank J. Bourgain, M. Combescure, J. Lebowitz, T. Spencer and M. Weinstein
for useful conversations. W.-M. Wang thanks Rutgers University, where
part of this work was done, for its hospitality. This work is partially supported 
by Rutgers center for non-linear analysis and NSF grants DMS-0100490 and DMR-9813268.
\endthanks 
\abstract 
We prove that at large disorder, Anderson localization in $\Z^d$ is stable under
localized time-periodic perturbations by proving that the associated quasi-energy
operator has pure point spectrum. The formulation of this problem is motivated by
questions of Anderson localization for non-linear Schr\"odinger equations.
\endabstract
\endtopmatter
\vfill \eject
\heading I. Introduction \endheading
Anderson localization for time {\it independent} random Schr\"odinger operators at
large disorder has been well known since the seminal work of Fr\"ohlich-Spencer
\cite{FS}. It is a topic with an extensive literature \cite{GMP, FMSS, vDK, AM, AFHS, AENSS},
to name a few.

Time-independent random Schr\"odinger operator is an operator of the form
$$H_0=\Delta+\gamma V,\,\tag 1.1$$
on $L^2(\R^d)$ or $\ell^2({\Z}^d)$, where $\Delta$ is the continuum or discrete
Laplacian, $\gamma$ is a positive parameter and $V$ is a random potential. We specialize
to discrete random Schr\"odinger operator. $H_0$ is then defined as the operator:
$$H_0=\Delta+\gamma V,\, \text {on}\, \ell^2({\Z}^d),\tag 1.2$$
\noindent
where the matrix element $\Delta_{ij}$, for $i$, $j\in\Z^d$ verify
$$\aligned \Delta_{ij}&=1\quad |i-j|_{\ell^1}=1\\
&=0\quad \text{otherwise};\endaligned\tag 1.3$$
\noindent
$\gamma$ is a positive parameter, the potential function $V$ is a diagonal matrix:
$V=\text{diag}(v_j),\, j\in{\Z}^d$, where $\{v_j\}$ is a family of independently
identically distributed (iid) real random variables with distribution $g$.
From now on, we write $|\,|$ for the $\ell^1$ norm: $|\,|_{\ell^1}$ on $\Z^d$. 
We denote $\ell^2$ norms by $\Vert\,\Vert$. 
The probability space $\Omega$ is taken to be $\R^{\Z^d}$ and the measure $P$ is 
$\prod_{j\in {\Z}^d}g(dv_j)$.

Let $\sigma (H_0)$ denote the spectrum of $H$. The spectrum can be decomposed 
into $\sigma_{\text{pp}}(H)$, $\sigma_{\text{ac}}(H)$ and
$\sigma_{\text{sc}}(H)$, where $\sigma_{\text{pp}}(H)$, the``pure point" spectrum, denotes the {\it closure} of
the set $S(H)=\{\lambda|\lambda\, \text{ is an eigenvalue of}\,H\}$,  $\sigma_{\text{ac}}(H)$ the absolutely
continuous spectrum and $\sigma_{\text{sc}}(H)$ the singular continous spectrum. We have the 
well established fact that $\sigma(H)$ and its decompositions $\sigma_{\text{pp}}$, $\sigma_{\text{ac}}$ and
$\sigma_{\text{sc}}$ are almost surely constant sets in $\R$, (see e.g., \cite{CFKS,PF}). 

\noindent {\it Remark.} The definition of $\sigma_{\text{pp}}$ is different from the usual definition, e.g.,
from that in \cite{RS}. This is in order that we have the stability property regarding spectral
decomposition as mentioned above. 

As is well known, $\sigma(\Delta)=[-2d, 2d]$. Let $\supp\, g$ be the support of
$g$, we know further (see e.g., \cite{CFKS,PF}) that 
$$\sigma(H)=[-2d, 2d]+\gamma\supp\, g\quad a.s.\tag 1.4$$
\ssk
The basic result proven in the references mentioned earlier is that under certain regularity
conditions on $g$, for $\gamma>>1$, and in any dimension $d$, the spectrum of $H_0$ is
almost surely pure point with exponentially localized eigenfunctions. This is called
{\it Anderson localization}, after the physicist P. W. Andreson \cite{An}. Physically this
corresponds to a lack of conductivity due to the localization of electrons. Anderson was the
first one to explain this phenomenon on theoretical physics ground.

The study of electron conduction is a many body problem. One needs to take into account
the interactions among electrons. This is a hard problem. The operator $H_0$ defined
in (1.2) corresponds to the so called 1-body approximation, where the interaction
is approximated by the potential $V$. The equation governing the system is
$$i\frac{\partial}{\partial t}\psi=(\Delta+\gamma V)\psi\tag 1.5$$
on $\Z^d\times[0,\infty)$.

This is the usual Schr\"odinger equation with a random potential. Since $V$ is independent
of $t$, the study of (1.5) could be reduced to the study of spectral properties of $H_0$.
Hence the importance of spectral results on $H_0$ mentioned earlier.

In this paper, we consider (1.5) perturbed by a bounded, localized (in space), time-periodic
potential. We study the equation:
$$i\frac{\partial}{\partial t}\psi=(\Delta+\gamma V+\lambda\Cal W)\psi\tag 1.6$$
on $\Z^d\times [0,\infty)$, where $V$ is as in (1.2), $\{v_j\}$ is a family of
(time-independent) i.i.d. random variables; $\Cal W=\Cal W(t,j)$, which we further
assume to be of the form: $$\Cal W(t,j)=\cos 2\pi(\omega t+\theta) W(j)\,(\omega>0).\tag 1.7$$
The motivation for studying (1.6) comes from questions of Anderson localization for non-linear
Schr\"odinger equations (see e.g., \cite {DS, FSW}), which in turn is a approximation to the
many body problem of electron conduction mentioned earlier. 

To proceed further, we assume
\item {(H1)} $g$ has bounded support
\item {(H2)} $g$ is absolutely continuous with a bounded density $\tilde g$:
$$g(dv)=\tilde g(v) dv, \, \Vert \tilde g\Vert _{\infty}<\infty$$
\item {(H3)} $$|W(j)|\leq Ce^{-b(\log\gamma)|j|}\, (C>0, b>0)$$
\noindent{\it Remark.} To prove Theorem 1.1 below, we only need $W$ to decay exponentially away
from the origin. In (H3), the $\gamma$ dependence of the rate of decay of $W$ is chosen to
coincide with the $\gamma$ dependence of the rate of decay of eigenfunctions of $H_0$, as our
main motivation for studying (1.6) comes from non-linear Schr\"odinger. We assume (H1) for
convenience. In the case $g$ is unbounded, we believe that supplemented with Lifshitz tail
arguments (see e.g., \cite{PF}), the proof presented in sect. II and III would go through.

We further note that due to the presence of the parameters $\gamma$ and $\lambda$, without 
loss, (H1,3) could be replaced by
\item {(H1')} $\supp g\subset [-1,1]$ 
\item {(H3')} $|W(j)|\leq e^{-b(\log\gamma)|j|}\, (b>0)$

\noindent Assume (H2), the precise spectral property of $H_0$ mentioned earlier is the
following (see e.g., \cite{vDK}):
\proclaim {Localization Theorem} Let $I$ be an
interval in $\R$. There exists $m>0$, such that for sufficiently large $\gamma$,
with probability $1$:

$\bullet$ $\sigma (H_0)\cap I$ is pure point, 

$\bullet$ the eigenfunctions $\psi_E$ corresponding to eigenvalues $E$ in $I$ satisfy
$$\liminf_{\Vert i\Vert\to\infty}-\frac{\log|\psi_E(i)|}{\Vert i\Vert}\geq m\log\gamma.$$ 
\endproclaim 

Define $$H(\theta)=H_0+\lambda\Cal W(\theta),\tag 1.8$$
the operator that first appeared in the RHS of (1.6). We consider $H(\theta)$ as a family
of Hamiltonians depending parametrically on the initial point $\theta\in\Cal S^1$, a unit
circle. We define $T_t$ to be the shift operator: $$T_t\theta=\omega t+\theta,\tag 1.9$$
and $\theta(t)=T_t\theta\in\Cal S^1$.

Let $U(t,s;\theta)$ ($t$, $s\in\R$) be the corresponding propagator: if at time $s$, a
solution to (1.6) is $\psi(s)$, then 
$$\psi(t)=U(t,s;\theta)\psi(s)\tag 1.10$$
is the solution at time $t$. For (1.6), using well known arguments, see e.g.,
\cite {Ho1, Ya}, we know that $U(t,s;\theta)$ is unitary and strongly continuous in $t$ and
$s$. Moreover, it satisfies $$U(t+a, s+a;\theta)=U(t, s; T_a\theta),\, (a\in\R).\tag
1.11$$

As in the usual construction (see e.g., \cite{Ho1, Ya}), we consider the enlarged space
$$\Cal K=\ell^2(\Z^d)\otimes L^2(\Cal S^1),\tag 1.12$$
and the one-parameter family of operators $\tilde U(t)$ ($t\in\R$) acting on $\Psi\in\Cal K$
by 
$$\aligned [\tilde U(t)\Psi](\theta)&=U(0,-t;\theta)[\Cal T_{-t}\Psi](\theta)\\
&=\Cal T_{-t}U(t,0;\theta)\Psi(\theta)\endaligned\tag 1.13$$
with 
$$[\Cal T_{-t}\Psi](\theta)=\Psi(T_{-t}\theta).\tag 1.14$$ 

It can be shown that the $\tilde U(t)$ here is a strongly continuous family of unitary
operators, see e.g., \cite{Ya}. By Stone's theorem, it can therefore be represented as
$$\tilde U(t)=e^{-iKt},\tag 1.15$$
where
$$K=\frac{\omega}{i}\frac{\partial}{\partial\theta}+\Delta+\gamma V+\lambda\cos 2\pi\theta
W\tag 1.16$$
on $\ell^2(\Z^d)\times L^2(\Cal S^1)$ is the quasi-energy operator. When $t=T=1/\omega$,
the period of the system, $U(T,0;0)$ is the Floquet operator. Formally, the generalized
eigenvalues and eigenfunctions of $U(T,0;0)$ and $K$ are related by (see e.g., \cite {JL})
$$\aligned K\psi&=\lambda\psi\\
U(T,0;0)\phi&=e^{-i\lambda T}\phi\\
\psi(\theta)&=e^{i\lambda \theta}U(\theta,0;0)\phi(0)\endaligned\tag 1.17$$

Our main result is
\proclaim{Theorem 1.1} Assume $g$ satisfies (H1',2) and $W$ satisfies (H3'). Let $I$ be an
interval in $\R$. There exists $a>0$, such that for sufficiently large $\gamma$,
with probability $1$:

$\bullet$ $\sigma (K)\cap I$ is pure point, 

$\bullet$ the eigenfunctions $\psi_E$ corresponding to eigenvalues $E$ in $I$ satisfy
$$\liminf_{\Vert i\Vert\to\infty}-\frac{\log|\psi_E(i)|}{\Vert i\Vert}\geq a\log\gamma.\tag
1.18$$ for all $\theta\in\Cal S^1$.
\endproclaim
\noindent{\it Remark.} Theorem 1.1 is deduced from localization properties of $H_0$. For
related KAM type of method used to study perturbations of dense pure point spectrum, see e.g.,
\cite {Ho2, TW}.

Using (1.13-1.15), Theoorem 1.1 and that any $\phi\in\ell^2(\Z^d)$ can be embedded
in $\Cal K=\ell^2(\Z^d)\otimes L^2(\Cal S^1)$ as $\psi=\phi\otimes 1$, we obtain
\proclaim{Corollary 1.2} Assume $g$ satisfies (H1',2) and $W$ satisfies (H3'). For sufficiently
large $\gamma$, for all $\phi\in\ell^2(\Z^d)$, all $\epsilon>0$, there exists $R>0$,
such that
$$\sup_t\sum_{|j|>R}|(U(t,0;\theta)\phi)(j)|^2<\epsilon\quad a.s.\tag 1.19$$
\endproclaim

Corollary 1.2 implies that Anderson localization is {\it stable} under bounded, localized
time periodic perturbations. It is a type of quantum stability result. For other 
results of related interests, see e.g., \cite{Be, Co, Sa}.  

Finally, we sketch the main ideas to prove Theorem 1.1. As in all other proofs of localization,
we use the established mechanism. We follow most closely \cite{vDK}. The operator
$K_0=K(\lambda=0)$ plays an important role, as in Fourier space:
$$K_0=2\pi n\omega+\Delta+\gamma V,\, n=0,\pm 1,\pm 2...\tag 1.20$$
We need two ingredients, both probabilistic in nature:
\item{(i)} A wegner estimate on regularity of eigenvalue spacing;
\item{(ii)} An initial localization estimate on finite volume Green's function.

In sect. II, we prove (i), which can be reduced to an estimate on the number of eigenvalues 
of $K_{\Lambda}$ ($K$ restrcted to a finite set $\Lambda\subset\Z^d$) in a given spectral 
interval. The bound is obtained by using the Helffer-Sj\"ostrand representation of
$f(K_{\Lambda})$ for $f\in C_0^{\infty}(\R)$.

In sect. III, we prove (ii). The initial localization estimate on $(E-K_{\Lambda})^{-1}$
is obtained from initial localization estimate on $(E-K_{0,\Lambda})^{-1}$ and localization
properties of $W$ in (H3'). The initial estimate on $(E-K_{0,\Lambda})^{-1}$ is provided
by {\it uniform} localization estimates on $H_0$. The fact that $\sigma(\Delta)$ is
bounded for the discrete Laplacian plays an essential role here. In both proofs of (i) and
(ii), we use Hilbert-Schmidt properties of $(E-K_{\Lambda})^{-1}$ and $(E-K_{0,\Lambda})^{-1}$.
\vskip 1cm
\heading II. Wegner Estimate for $K$ \endheading
Recall from sect. I, the quasi-energy operator $K$:
$$\aligned K&=\frac{\omega}{i}\frac{\partial}{\partial\theta}+H(\theta)\\
&=\frac{\omega}{i}\frac{\partial}{\partial\theta}+\Delta+\gamma V+\lambda\cos
2\pi\theta W\endaligned $$
on $\ell^2(\Z^d)\times L^2(\Cal S^1)$, as defined in (1.16).

To prove Theorem 1.1, we proceed in the usual way. We need a Wegner estimate on regularity
of eigenvalue spacing and an initial estimate on the Green's function. Toward that end,
let $\Lambda$ be a finite subset in $\Z^d$. Define
$$\aligned \Delta_{\Lambda}(i,j)&=\Delta(i,j)\quad \text {if}\, i,\,j\in\Lambda\\
&=0\quad \text{otherwise};\endaligned$$
and $$ K_{\Lambda}=\frac{\omega}{i}\frac{\partial}{\partial\theta}+\Delta_\Lambda+\gamma
V+\lambda\cos 2\pi\theta W$$
on $\ell^2(\Lambda)\times L^2(\Cal S^1)$.

In this section, we prove the Wegner estimate. By using the standard shift (in V) argument, see
e.g., the proof of Proposition 3.1 in \cite{W}, we know that for $\epsilon<<1$
\proclaim {Proposition 2.1}
$$\text {Prob}(\text{dist}(E,\sigma(K_{\Lambda})\leq\epsilon)\leq
C\frac{N_I(E)\epsilon|\Lambda|}{\gamma}\tag 2.1$$
where $I=(E-1,E+1)$, $N_I(E)$ is the number of eigenvalues in $I$.
\endproclaim
Applying Proposition 2.1 to $K_{\Lambda}$, we have
\proclaim {Lemma 2.2 (Wegner estimate for $K_\Lambda$)}
$$\text {Prob}(\text{dist}(E,\sigma(K_{\Lambda})\leq\epsilon)\leq
C\big (\frac{\epsilon}{\omega}\big )|\Lambda|^4 \big (\frac{2d+\gamma}{\gamma}\big)\tag 2.2$$
where $C$ only depends on $\lambda$ and the probability distribution $g$.
\endproclaim
\demo{Proof} Let $f\in C_0^{\infty}(\R;\R^+)$, $f(x)=1$, $x\in I$, $f(x)=0$ if
$x\geq E+2$ or $x\leq E-2$. Then
$$\aligned N_I(E)&\leq \Tr f(K_\Lambda)\\
&=\Tr \{\big ( \frac{i}{2\pi}\big)\int \partial_{\bar z}\tilde f(z)(z-K_\Lambda)^{-1} d\bar
z\wedge dz\}\endaligned\tag 2.3$$
where $\tilde f \in C_0^{\infty}(\C)$ is an almost analytic extension of $f$, i.e.,
$\tilde f=f$ on $\R$ and $\partial_{\bar z}\tilde f$ vanishes on $\R$ to infinite order, 
see e.g., \cite{HS, D}. Let
$$\aligned \Cal W&=\lambda\cos 2\pi\theta W\\
K_{0,\Lambda}&=K_\Lambda-\Cal
W=\frac{\omega}{i}\frac{\partial}{\partial\theta}+\Delta_\Lambda+\gamma V.\endaligned\tag 2.4$$

For simplicity of notation, we write $K_0$ for $K_{0,\Lambda}$, $K$ for $K_\Lambda$. We
note that $\frac{\omega}{i}\frac{\partial}{\partial\theta}$ commutes with
$H_\Lambda=\Delta_\Lambda+\gamma V$. Passing to the dual variable of $\theta$ by Fouries
series, (and abusing the notation), we have 
$$\aligned \Cal W&=\lambda\big (\frac {T_++T_-}{2}) W\\
K_0&=2\pi n\omega+\Delta_\Lambda+\gamma V\quad n=0,\pm 1,\pm 2...\endaligned\tag 2.5$$
where $T_{\pm}$ are unit shift operators on $\Z$: 
$$(T_{\pm}f)(n)=f(n\pm 1).\tag 2.6$$

Using the resolvent equation twice, we have
$$\aligned (z-K)^{-1}=(z&-K_0)^{-1}-(z-K_0)^{-1}\Cal W (z-K_0)^{-1}\\
&+(z-K_0)^{-1}\Cal W (z-K_0)^{-1}\Cal W (z-K)^{-1}.\endaligned\tag 2.7$$
Since $(z-K_0)^{-1}$ is diagonal in $n$, $\Cal W$ only has off-diagonal elements, the second
term in the RHS of (2.7) is traceless. Substituting (2.7) into (2.3), we then obtain
$$\aligned N_I(E)\leq 
&\Tr \{\big ( \frac{i}{2\pi}\big)\int \partial_{\bar z}\tilde f(z)(z-K_0)^{-1} d\bar
z\wedge dz\}\\
&+\Tr \{\big ( \frac{i}{2\pi}\big)\int \partial_{\bar z}\tilde f(z)(z-K_0)^{-1}\Cal W
(z-K_0)^{-1}\Cal W (z-K)^{-1} d\bar z\wedge dz\}\\
\overset\text{def}\to=&I_1+I_2 \endaligned\tag 2.8$$

We first evaluate $I_1$: 
Since $\sigma(\Delta)\subset[-2d,2d]$, we have 
$$\sigma(\Delta+\gamma V)\subset[-2d-\gamma,2d+\gamma].\tag 2.9$$
Recall that $\supp f=(E-2,E+2)$. So we can take $\supp\Re \tilde f=(E-2,E+2)$; i.e., 
$\Re z\in (E-2,E+2)$.

When evaluating the trace in $I_1$, we only need to sum over $n$, such that 
$$\aligned |z-2\pi n\omega|&\leq 2d+\gamma+1\\
|n-\frac{z}{2\pi\omega}|&\leq \frac{1}{2\pi\omega}(2d+\gamma+1)\\
|n-\frac{E}{2\pi\omega}|&\leq \frac{1}{2\pi\omega}(2d+\gamma+3).\endaligned\tag 2.10$$
Otherwise $z-K_0$ is invertible, the integrand is analytic in $z$ and the integral is
$0$ for such $n$ by using Stokes' formula. Hence
$$\aligned I_1=\big ( \frac{i}{2\pi}\big)\sum_{j\in\Lambda}&\sum_{|n-\frac{E}{2\pi\omega}|
\leq\frac{1}{2\pi\omega}(2d+\gamma+3)}\\
&\int \partial_{\bar z}\tilde f(z)(z-2\pi
n\omega-\Delta-\gamma V)^{-1}(j,j) d\bar z\wedge dz.\endaligned\tag 2.11$$

As is the standard practice, we split the $d\bar z\wedge dz$ integration into
$|\Im z|\geq\alpha$ and $|\Im z|\leq\alpha$ for some $\alpha>0$ to be chosen conveniently.
So
$$\aligned |I_1|\leq &\big (\frac{1}{2\pi}\big)|\Lambda|\frac{2(2d+\gamma+3)+1}{2\pi\omega}\\
&\sup_{j\in\Lambda}\big (|\int_{|\Im z|\geq\alpha} \partial_{\bar z}\tilde f(z)(z-2\pi
n\omega-\Delta-\gamma V)^{-1}(j,j) d\bar z\wedge dz|\\
&\qquad +|\int_{|\Im z|\leq\alpha} \partial_{\bar z}\tilde f(z)(z-2\pi n\omega-\Delta-\gamma
V)^{-1}(j,j) d\bar z\wedge dz|\big )\\
=&\frac{|\Lambda|}{4\pi^2\omega}(2(2d+\gamma+3)+1)\big(\Cal O(1)\frac{1}{\alpha}+
\Cal O_M(1)|\Im z|^M\frac{1}{|\Im z|}\big )\endaligned\tag 2.12$$
for all $M\in\N^+$, where we used $|\partial_{\bar z}\tilde f(z)|\leq \Cal O_M(1)|\Im z|^M$
for all $M$ and self-adjointness. 

Choosing $M=\alpha=1$, we then obtain
$$|I_1|\leq \Cal O(1)\frac{|\Lambda|}{\omega}(2d+\gamma),\tag 2.13$$
where $\Cal O(1)$ is uniform in $E$, $\omega$, $d$ and $\gamma$.

We now estimate $I_2$, where the main complication comes from the term $(z-K)^{-1}$. The
only control we have is via $\Im z$. We split the sum over $n$ similar to (2.11). Anticipating
ahead, we split the sum into $|n-\frac{E}{2\pi\omega}|\leq\frac{1}{2\pi\omega}(2d+\gamma+3)+1$
and its complement:
$$\aligned &I_2\\
=&\big (\frac{i}{2\pi}\big)\sum_{j\in\Lambda}
\big (\sum_{|n-\frac{E}{2\pi\omega}|\leq\frac{1}{2\pi\omega}(2d+\gamma+3)+1}\int
\partial_{\bar z}\tilde f(z)[(z-K_0)^{-1}\Cal W (z-K_0)^{-1}\Cal W (z-K)^{-1}](n,j,n,j) d\bar
z\wedge dz\\ 
&\qquad
+\sum_{|n-\frac{E}{2\pi\omega}|>\frac{1}{2\pi\omega}(2d+\gamma+3)+1}\int \partial_{\bar
z}\tilde f(z)[(z-K_0)^{-1}\Cal W (z-K_0)^{-1}\Cal W (z-K)^{-1}](n,j,n,j) d\bar z\wedge dz\big
)\\ =&r_1+r_2.\endaligned\tag 2.14$$

$r_1$ can be estimated in the same way as in (2.12). We write out the kernel:
$$\aligned &[(z-K_0)^{-1}\Cal W (z-K_0)^{-1}\Cal W (z-K)^{-1}](n,j,n,j)\\
=&\sum_{n',n'',j',j''}(z-2\pi n\omega-\Delta-\gamma V)^{-1}(n,j,n,j')\Cal W(n,j',n',j')\\
&\qquad (z-2\pi n'\omega-\Delta-\gamma V)^{-1}(n',j',n',j'')\Cal W(n',j'',n'',j'')\\
&\qquad (z-K)^{-1}(n'',j'',n,j).\endaligned\tag 2.15$$
We note that from (2.5), $|n-n'|=1$, $|n''-n|\leq 2$. Taking $M=3$ instead of $1$ and
summing over $j$, $j'$, $j''$, $n$, $n'$, $n''$, we obtain
$$|r_1|\leq\Cal O(1)\frac{|\Lambda|^3}{\omega}(2d+\gamma).\tag 2.16$$

Estimation of $r_2$ is different from that of $I_1$, as a priori we cannot conclude that the
integrand is analytic in $z$ for such large $n$. Instead, we do the following:
$$\aligned |r_2|\leq\Cal
O(1)\sum_{j,j',j''}&\sum_{|n-\frac{E}{2\pi\omega}|\geq\frac{1}{2\pi\omega}(2d+\gamma+3)+1}\\
&\int|\partial_{\bar z}\tilde f(z)|\sum_{n',n''}|[(z-K_0)^{-1}\Cal W (z-K_0)^{-1}\Cal
W](n,j,n'',j'')|\frac{1}{|\Im z|}d\bar z\wedge dz.\endaligned\tag 2.17$$
Using the fact that $|n-n'|=1$, $|n''-n|\leq 2$, the sum over $n$, $n'$, $n''$ is convergent.
We obtain
$$|r_2|\leq\Cal O(1)|\Lambda|^3.\tag 2.18$$

Combining (2.16, 2.18) with (2.13) in (2.8), we have 
$$N_I(E)\leq\Cal O(1)\frac{|\Lambda|^3}{\omega}(2d+\gamma),\tag 2.19$$
where $\Cal O(1)$ is uniform in $E$, $d$, $\Lambda$, $\omega$ and $\gamma$.
Substituting (2.19) into (2.1), we obtain the lemma.
\hfill {$\square$}
\enddemo
\vskip 1cm

\heading III. Initial estimate for localization and proof of Theorem 1.1\endheading
The initial estimate for localization for $K$ is deduced from localization estimates on
$H_0=\Delta+\gamma V$:
\proclaim {Proposition 3.1} There exist $a>0$, $\gamma_0>0$,
$L_0>0$, such that if we let
$L_{n+1}=L_n^\alpha$ ($1<\alpha<2$), $i\in\Z^d$, $\Lambda_n=[-L_n,L_n]^d+i$, then
for $\gamma>\gamma_0$, any set $S\subset [-2d-\gamma, 2d+\gamma]$ of
$\Cal O(1)\gamma$ elements, with probability $\geq 1-\frac{1}{L_n^p}$ ($p>2d$), for all
$j_n\in\partial\Lambda_n$, all $E\in S$,
$$|(E-H_{\Lambda_n})^{-1}(i,j_n)|\leq Ce^{-a\log\gamma|i-j_n|}\quad (C>0,\, a>0).\tag 3.1$$
\endproclaim
\demo {Proof} (3.1) is obtained by patching together the usual localization proof, see e.g.,
\cite{FMSS}. We will thus only mention that aspect. As in all large disorder case, $L_0=\Cal
O(1)$. The Wegner estimate for $H_0$ is:
$$\text {Prob}(\text{dist}(E,\sigma(H_{\Lambda})\leq\epsilon)\leq
\frac{C\epsilon|\Lambda|}{\gamma},\tag 3.2$$
where $C$ only depends on the distribution $g$, see e.g., \cite{vDK}. To get (3.1)
for $L_0$, we take $$\epsilon=\Cal O(1)\gamma^q\quad (0<q<1/2).\tag 3.3$$ So using (3.2), we have 
$$\Vert (E-H_{\Lambda_{L_0}})^{-1}\Vert\leq \Cal O(1)\gamma^{-q}\quad (0<q<1/2)\tag 3.4$$
with probability $\geq 1-\Cal O(1)\gamma^{-1+q}.$ 
\noindent(Recall that $L_0=\Cal O(1)$.) 

Let $\Lambda_{L_0(i)}$, $\Lambda_{L_0(j)}$ be 2 cubes of side length $2L_0$, centered at
$i$, $j\in\Z^d$. We note that if $\Lambda_{L_0}(i)\cap\Lambda_{L_0}(j)=\emptyset$, then
$$\aligned&\text{Prob}\{\text{both}\Vert (E-H_{\Lambda_{L_0}(i)})^{-1}\Vert\geq \Cal O(1)\gamma^{-q}
\,\text{and} \Vert (E-H_{\Lambda_{L_0}(i)})^{-1}\Vert\geq \Cal O(1)\gamma^{-q}\}\\
&\qquad\leq\Cal O(1)\gamma^{-2(1+q)}<<\gamma^{-1}\endaligned\tag 3.5$$
for $0<q<1/2$.

Let $\tilde L_0=\gamma^\delta$, $\delta>0$ only depends on $d$, $q$ and $p$. Using (3.2, 3.5),
we obtain that $$|(E-H_{\Lambda_{\tilde L_0(i)}})^{-1}(i,j)|\leq e^{-a'\log\gamma|i-j|}\quad
(a'>0),\tag 3.6$$
for all $j\in\partial\Lambda_{\tilde L_0(i)}$ with probability $$\geq 1-\Cal O(1)\tilde
L_0^{-p}\quad (p>2d).$$ 
for any set $S\subset[-2d-\gamma, 2d+\gamma]$ with $\Cal O(1)\gamma$ elements. 

Renaming $\tilde L_0$ as the new $L_0$,
using (3.2, 3.6) as our initial input in the localization mechanism we obtain that
$$|(E-H_{\Lambda_n})^{-1}(i,j_n)|\leq Ce^{-a\log\gamma|i-j_n|}\quad (C>0,\, a>0)\tag 3.7$$
with probability $\geq 1-L_n^{-p}$, ($p>2d$), for all $E\in S\subset [-2d-\gamma, 2d+\gamma]$. 
\hfill {$\square$}
\enddemo

Let $K_{0,\Lambda}$ be defined as in (2.4). We have 
\proclaim {Lemma 3.2}  There exist $a>0$, $L>0$, such that if $\gamma>>1$ and if we let 
$i\in\Z^d$, $\Lambda=[-L,L]^d+i$, then for all $j\in\partial\Lambda$, $x$, $y\in [0,1)$
$$|(E-K_{0,\Lambda})^{-1}(i,x;j,y)|\leq C
\big(\frac{\gamma+2d}{\omega}\big)e^{-a\log\gamma|i-j|}\quad (C>0,\, a>0)\tag 3.8$$
with probability $\geq 1-1/L^p$ ($p>2d$).
\endproclaim
\demo{Proof} Using Fourier series, we have
$$(E-K_{0,\Lambda})^{-1}(i,x;j,y)=\sum_{n=0,\pm 1,...}(E-2\pi
n\omega-H_0)^{-1}(i,j,n)e^{in(x-y)}.\tag 3.9$$
$$\aligned |(E-K_{0,\Lambda})^{-1}(i,j;x-y)|&\leq \sum_{n=0,\pm 1,...}|(E-2\pi
n\omega-H_0)^{-1}(i,j,n)|\\
&\leq \sum_{|n-\frac{E}{2\pi\omega}|\leq\frac{1}{\pi\omega}(2d+\gamma)}
|(E-2\pi n\omega-H_0)^{-1}(i,j,n)|\\
&\quad + \sum_{|n-\frac{E}{2\pi\omega}|>\frac{1}{\pi\omega}(2d+\gamma)}
|(E-2\pi n\omega-H_0)^{-1}(i,j,n)|\\
&\leq \Cal O(1)\frac{2d+\gamma}{\omega} e^{-a\log\gamma|i-j|}+\Cal
O(1)\frac{1}{\omega} e^{-a\log(2d+\gamma)|i-j|}\\
&\leq C \frac{2d+\gamma}{\omega} e^{-a\log\gamma|i-j|}\endaligned\tag 3.10$$
with probability $\geq 1-1/L^p$, where we assumed $\gamma>>1$ and used (3.1) in estimating the
first sum and standard elliptic estimate on the second sum.
\hfill {$\square$}
\enddemo

In order to prove Proposition 3.4 below, we need a slight generalization of Lemma 3.2,
which we state without proof as 
\proclaim {Corollary 3.3}
There exist $a>0$, $L>0$, such that if $\gamma>>1$ and if we let 
$i\in\Z^d$, $\Lambda=[-L,L]^d+i$, then for all $j,\,j'\in\Lambda$, $|j-j'|\geq L/4$,
$y$, $y'\in [0,1)$
$$|(E-K_{0,\Lambda})^{-1}(i,x;j,y)|\leq C
\big(\frac{\gamma+2d}{\omega}\big)e^{-a\log\gamma|i-j|}\quad (C>0,\, a>0)\tag 3.11$$
with probability $\geq 1-1/L^p$ ($p>2d$), the $p$ here is not necessarily the same as in 
Lemma 3.2.
\endproclaim

Using assumption (H3') and Proposition 3.1, we are now ready to prove
\proclaim {Proposition 3.4 (Initial estimate for $K_\Lambda$)} 
There exist $\tilde a>0$, $L\in\N^+$, such that if we let $i\in\Z^d$, and
$\Lambda=[-L,L|^d+i$, then for $\gamma>>1$, all $j\in\partial\Lambda$, $x$, $y\in[0,1)$
$$|(E-K_{\Lambda})^{-1}(i,x;j,y)|\leq
C\frac{\gamma+2d}{\omega}e^{-\tilde a\log\gamma|i-j|}\quad (C>0,\, \tilde a>0).\tag 3.12$$ with
probability $\geq 1-\frac{1}{L^p}$ ($p>2d$)
\endproclaim
\demo{Proof} We deduce (3.12) from (3.11) by using the resolvent equation, the Wegner estimate
in (2.2) and localization property of $W$ in (H3'). Iterating the resolvent equation twice,
we have (writing $K$ for $K_\Lambda$, $K_0$ for $K_{0,\Lambda}$):
$$\aligned (E-K)^{-1}(i,x;&j,y)=(E-K_0)^{-1}(i,x;j,y)\\
-&[(E-K_0)^{-1}\Cal W (E-K_0)^{-1}](i,x;j,y)\\
+&[(E-K_0)^{-1}\Cal W (E-K)^{-1}\Cal W (E-K_0)^{-1}](i,x;j,y)\\
\overset\text{def}\to=&I_1+I_2+I_3.\endaligned\tag 3.13$$

We use (3.8) to estimate the first term in the RHS of (3.13). To estimate the second term,
we use (3.8, 3.11). Let $\tilde b=\epsilon\min(a,b)$ for some $\epsilon>0$ to be determined
later. We write
$$\aligned I_2=&e^{-\tilde b\log\gamma|i-j|}I_2e^{\tilde b\log\gamma|i-j|}\\
\overset\text{def}\to=&e^{-\tilde b\log\gamma|i-j|}\tilde I_2.\endaligned\tag 3.14$$
We only need to bound $\tilde I_2$.
$$\aligned |\tilde I_2|&\leq e^{\tilde b\log\gamma|i-j|}\sum_k\int dt|[(E-K_0)^{-1}|\Cal
W|^{1/2}] (i,x;k,t)|\\
&\qquad |[|\Cal W|^{1/2}(E-K_0)^{-1}] (k,t;j,y)|\\
&\leq e^{\tilde b\log\gamma|i-j|}\sum_k\big[\{\int dt
[(E-K_0)^{-1}(i,x;k,t)]^2|W(k)|\}^{1/2}\\ &\quad\{\int dt
[(E-K_0)^{-1}(k,t;j,y)]^2|W(k)|\}^{1/2}\big]\\ &\leq e^{\tilde
b\log\gamma|i-j|}\sum_k\big[\{\sum_{n=0,\pm 1,\pm 2...} [(E-2\pi
n\omega-H_0)^{-1}(i,k)]^2|W(k)|\}^{1/2}\\ &\quad\{\sum_{n'=0,\pm 1,\pm 2...} [(E-2\pi
n'\omega-H_0)^{-1}(k,j)]^2|W(k)|\}^{1/2}\big]\\ &\leq\Cal
O(1)\frac{2d+\gamma}{\omega}|\Lambda|^q\endaligned\tag 3.15$$ for some $q>0$, with probability
$\geq 1-\frac{1}{L^p}$ ($p>2d$), where we estimated the sum over $n$ similar to (3.10), and we
used (3.11) and the Wegner estimate for $H_0$ in (3.2) with $\epsilon=|\Lambda|^{-q}$, $q$
adjusted according to $p$.

We now estimate $I_3$. Similar to (3.14), we write
$$\aligned I_3=&e^{-\tilde b\log\gamma|i-j|}I_3e^{\tilde b\log\gamma|i-j|}\\
\overset\text{def}\to=&e^{-\tilde b\log\gamma|i-j|}\tilde I_3.\endaligned\tag 3.16$$
$$\aligned |\tilde I_3|\leq& e^{\tilde b\log\gamma|i-j|}\sum_{k,k'}\int dt\int
dt'|[(E-K_0)^{-1}|\Cal W|^{1/2}] (i,x;k,t)||\Cal
W|^{1/2}[(E-K)^{-1}|\Cal W|^{1/2}] (k,t;k',t')|\\
&|[|\Cal W|^{1/2}(E-K_0)^{-1}] (k',t';j,y)|\\
\leq &e^{\tilde b\log\gamma|i-j|}\sum_{k,k'}\big[ \{\int dt\int dt'
|W(k)|[(E-K)^{-1}(k,t;k',t')]^2|W(k')|\}^{1/2}\\
&\{\int dt[(E-K_0)^{-1}(i,x;k,t)]^2|W(k)|\}^{1/2}\\ 
&\{\int dt'[(E-K_0)^{-1}(k',t';j,y)]^2|W(k')|\}^{1/2}\big]\\
\leq& e^{\tilde b\log\gamma|i-j|}\{\sum_{k,k'}\int dt\int dt'
|W(k)|[(E-K)^{-1}(k,t;k',t')]^2|W(k')|\}^{1/2}\\
&\{\sum_{k} \int dt[(E-K_0)^{-1}(i,x;k,t)]^2|W(k)|\}^{1/2}
\{\sum_{k'} \int dt'[(E-K_0)^{-1}(k',t';j,y)]^2|W(k')|\}^{1/2}\\
\overset\text{def}\to=&S_1S_2S_3\endaligned\tag 3.17$$

$S_2$, $S_3$ can be similarly estimated as in (3.15).
$$S_1\leq \Cal O(1)\Vert (E-K)^{-1}\Vert_{\text {HS}},\tag 3.18$$
where $\Vert \,\Vert_{\text {HS}}$ denotes the Hilbert-Schmidt norm. From the resolvent
equation, we have
$$(E-K)^{-1}=(E-K_0)^{-1}-(E-K_0)^{-1}\Cal W (E-K)^{-1}.\tag 3.19$$
To estimate the H-S norm, we sum over $n$ similar to (3.10, 3.16). Using (3.2, 2.2), we obtain
$$\aligned \Vert (E-K)^{-1}\Vert_{\text {HS}}&\leq
\Vert (E-K_0)^{-1}\Vert_{\text {HS}}(1+\lambda \Vert (E-K)^{-1}\Vert_{L^2})\\
&\leq\Cal O(1)|\Lambda|^s\frac{2d+\gamma}{\omega}\lambda\endaligned\tag 3.20$$
for some $s>0$, with probability $\geq 1-1/L^p$ ($s$ depends on $p$). 

Combining the estimates on $I_1$, $I_2$ and $I_3$ in (3.8), (3.14-3.20), adjusting $q$, $s$ and
$L$, we obtain (3.12).
\hfill {$\square$}
\enddemo
\ssk
\demo{Proof of Theorem 1.1} Using Lemma 2.2 and Proposition 3.4, Theorem 1.1 follows via the
standard route of localization proofs and polynomial boundedness of generalized eigenfunctions
of $K$, see e.g., \cite{Si}. (See also \cite{vDK, FMSS}.)
\hfill {$\square$}
\enddemo

\vfill \eject

\end